\documentclass[preprint,review,11pt]{elsarticle}

\usepackage{lscape, dashrule, amsfonts, newunicodechar}
\usepackage[utf8]{inputenc}
\usepackage{hyperref}
\newunicodechar{✓}{\checkmark}
\usepackage{amsmath, verbatim, amssymb, latexsym}
\usepackage{multirow, graphics, graphicx, color,hyperref}
\usepackage{psfrag, float, fancybox, caption,fancyhdr}

\usepackage[colorinlistoftodos,bordercolor=orange,backgroundcolor=orange!20,linecolor=orange,textsize=scriptsize]{todonotes}
\usepackage[pagewise]{lineno}

\usepackage{colortbl}
\usepackage{comment}
\usepackage{algpseudocode}
\usepackage{epstopdf}
\usepackage{natbib}
\usepackage{verbatim}
\usepackage{multirow}
\usepackage{supertabular}
\usepackage[]{algorithm2e}
\usepackage{url}            % simple URL typesetting
\usepackage{booktabs}       % professional-quality tables
\usepackage{amsfonts}       % blackboard math symbols
\usepackage{nicefrac}       % compact symbols for 1/2, etc.
\usepackage{microtype}      % microtypography
\usepackage[toc,page]{appendix}
\usepackage{comment}
\graphicspath{{./../figures/},{./../visio/}}
\makeatletter
\def\@author#1{\g@addto@macro\elsauthors{\normalsize%
    \def\baselinestretch{1}%
    \upshape\authorsep#1\unskip\textsuperscript{%
      \ifx\@fnmark\@empty\else\unskip\sep\@fnmark\let\sep=,\fi
      \ifx\@corref\@empty\else\unskip\sep\@corref\let\sep=,\fi
      }%
    \def\authorsep{\unskip,\space}%
    \global\let\@fnmark\@empty
    \global\let\@corref\@empty  %% Added
    \global\let\sep\@empty}%
    \@eadauthor={#1}
}
\makeatother

\journal{TBD}

\begin{document}

\begin{frontmatter}

\title{A Stochastic Unequal Area Facility Layout Problem}

\author{Afshin OroojlooyJadid\corref{cor1}\fnref{a}}
\cortext[cor1]{Corresponding Author; E-mail: oroojlooy@lehigh.edu; Phone: (484)666-8370}
\author{Mohammad Firouz\fnref{b}}

\address[a]{Industrial and Systems Engineering Departement, Lehigh University, Bthlehem, PA 18015.}
\address[b]{Operations Research Laboratory, General Motors Research and Development, Warren, MI 48092.}

\begin{abstract}
\noindent In this study we introduce a new method to solve the Dynamics Facility Layout Problems (DFLPs). To represent each layout, we use the slicing tree method integrated with our proposed heuristic to obtain promising initial solutions. Then, we consider the case of adding new departments into the current layout with stochastic flows. We use simulation to model the complexity of stochastic nature of the problem. To improve the quality of the initial solution a genetic algorithm is joined with the simulation module. Finally, to demonstrate the performance of our method, we solve several cases existing in the literature to show the efficiency of our algorithm. 
\end{abstract}

\begin{keyword}
Facility Layout Problem \sep Genetic Algorithm (GA) \sep Simulation \sep Stochastic \sep Unequal Area
\end{keyword}

\end{frontmatter}

% \linenumbers


\newfloat{display}{thp}{lod}
\floatname{display}{Display}

\section{Introduction}
\label{sec:intro}
In response to considerably fast technological changes, intensified competition, and shortened product life cycles, production systems need to increase their flexibility to change the number or the size of specific machines/departments in their layouts. 

Within the layout of a manufacturing system, flows of material is common between any two department. Such material flows are commonly represented by a from-to matrix. Once a new department is added or an existing department is modified, such material flows inevitably change throughout the whole layout with different degrees. Therefore the new material flows, are a function of the current material flows and a stochastic portion due to such change. 

Additionally, such changes happen at indefinite points of time. The time for such changes can vary significantly and is motivated by several reasons from maintaining market share in the current products, or competition with other manufacturers due to introduction of a new product, to even firm's strategic plans for changes in the long run. As such time of change is variable, the cost of change varies as well.

Traditional facility layout problems consider all the parameters including number of departments, material flow between them, and the related cost factors as constant, which is generally referred to as the Facility Layout Problem (FLP). The minimization objective function in FLP is most often the material flow cost. In a strategic view towards the FLPs, Dynamic Facility Layout Problem (DFLP) has emerged that considers a long term plan and has some rearrangement times at each period of production. At each period, there is a trade-off between the potential change in the flows between departments, and therefore the chance of getting a better layout, and the rearranging cost. We refer to \citet{kulturel2007approaches} for a complete review of the DFLPs. 

As mentioned, most firms have to rapidly adapt to the competitive market in order to maintain their marketshare. To do this, they need to make changes in the current products or add new ones as quick as possible which inevitably will lead to changes in the layout. Dealing with this condition, the new generation of DFLPs has emerged, that considers changes in the type and the number of facilities as well as changes in the material flow between them. First contribution to this type of DFLPs belongs to \citet{dong2009shortest}, with consideration of some new machines being added to the layout and some old machines removed from it in a deterministic setting.

In this paper, we study an Unequal Area Dynamic Facility Layout Problem UA-DFLP, considering departments' adding/removing under stochastic conditions and stochastic material flow changes among the departments. Such situation is observed often times in food and electronic industries. In these settings, the demand forecast for products can only be obtained as a lower and upper bound as [\emph{Lower-bound , Upper-bound}]. However, usually the center of this interval is used as product's demand forecast. Also, the range of the interval increases with higher inaccuracy comes from the market demand research so that even the center value will not remain stable. Therefore, using mean demand in deterministic DFLPs will encounter shortcomings from this perspective. %Additionally, the size of the departments that are added to the system are not necessarily equal which is another drawback to the current studies of DFLPs \citet{XXXX}. \todo{add reference.}

While there is a considerable contribution to DFLPs, up to our knowledge, no research has been done for our problem setting with unequal area departments added or removed in stochastic conditions and stochastic material flow and stochastic time of change. In this research we present an effective method for solving problems in such class of problems in DFLPs. 

The remainder of this paper is organized as follows. In section \ref{sec:litReview} we review the literature of the facility layout problem most closely related to our work. In section \ref{sec: model}, we introduce our model. Section \ref{sec: solution} describes our hybrid algorithm and proposes a novel methodology for creating promising initial solutions. To demonstrate the efficiency of the proposed algorithm we provide our results in section \ref{sec: computational}. Finally, section \ref{sec: conclusions} concludes the paper.

\section{Relevant Literature}
\label{sec:litReview}

The first step towards Unequal Area Facility Layout Problems (UA-FLPs) is taken by \citet{armour1963heuristic}. They consider departments with different areas and shapes, using the center of gravity of each department for calculating the distances without considering certain directions for departments. \cite{armour1963heuristic} minimize the transportation cost between departments with rectilinear distances. \citet{montreuil1991modelling} presents the first MIP-FLP model for the continuous-representation-based FLP in which they use a distance-based objective function similar to \citet{armour1963heuristic}. The presented model in \citet{montreuil1991modelling} is frequently referred to as the FLP0. 

\citet{meller1996facility} classifies and reviews previous contributions and redefines the problem, considering material handling cost as the objective function and constraints as followings:

\begin{enumerate}[(a)]
\item All the departments with certain area requirements should be located within the predetermined limited space.
\item There should not be any interference between departments in the final layout.
\item For each department there is a maximum ratio of length to width which should not be violated.
\end{enumerate}

The minimization objective function for linear distances in \citet{meller1996facility} is as follows:
\[\sum\limits_{j=1}^n\sum\limits_{i=1,i\neq j}^nf_{ij}d_{ij},\]
where $f_{ij}$ represents the material flow between departments $i$ and $j$ and $d_{ij}$ is the distance between them. All distances are measured center to center and are considered as rectilinear.

\citet{meller1998optimal} presents the first contribution to the UA-FLP using a Mixed Integer Programming (MIP) formulation by enhancing the old model presented by \citet{montreuil1991modelling}. They use an advanced substitute perimeter constraint to obtain optimal solution for problems with up to eight departments. However, their model still is not useful enough for large scale practical problems. \citet{sherali2003enhanced}, using a novel polyhedral outer approximation, enhanced the accuracy of the model presented by \citet{meller1998optimal} and were able to find optimal solution for problems with up to nine departments. \citet{logendran2006methodology} presents a mixed-binary nonlinear programming model with consideration of the geometry of the departments, using material handling costs as the objective function. \citet{logendran2006methodology} presents a tabu search algorithm to solve the problem and represents each layout using the slicing tree method \cite{otten1982automatic}. \citet{liu2007sequence} uses sequence-pair representation and presents a new formulation for UA-FLPs. \citet{liu2007sequence} are able to solve problems with up to to eleven departments.

\citet{anjos2006new} presents a two-phase mathematical programming formulation to the UA-FLPs. They use a promising starting point solution in the first phase to feed an exact formulation of the problem with a non-convex mathematical program with equilibrium constraints. In the second phase, an interactive algorithm solves the precise formulation of the problem. Table \ref{tab:exact} summarize the MIP models.

\begin{table}[]
\centering
\caption{Summary of the studies with exact solutions}
\label{tab:exact}
\begin{tabular}{lllc}
\hline
Paper               & Model                              & Year & Max Departments \\
\hline
\citet{montreuil1991modelling}          & FLP0                               & 1990 & -                         \\
\citet{meller1996facility}     & First MIP model                    & 1996 & 6                         \\
\citet{sherali2003enhanced}    & Polyhedral approximation           & 2003 & 9                         \\
\citet{castillo2005optimization}    & MINLP-$\epsilon$ accurate scheme           & 2005 & 9                         \\
\citet{anjos2006new} & two-phase mathematical programming & 2006 & -                         \\
\citet{liu2007sequence}     & First Sequence pair in MIP         & 2007 & 9 \\                       
\hline
\end{tabular}
\end{table}

FLP is known as an NP-Hard problem \citep{drira2007facility} and exact solutions cannot be found in polynomial time so that several heuristic and meta-heuristic algorithms are proposed to solve it.  \citet{azadivar2000facility} represents each layout using slicing tree and introduces a Gnetic Algorithm (GA) to get sufficiently good solutions. They consider the case that each machine has its unique queuing and breakdown distribution and conclude that simulation is very effective to evaluate each layout. Some important meta-heuristic methods that have been used for solving DFLPs are shown in Table \ref{tab:heuristic}. \citet{drira2007facility} also presents a helpful survey.

\begin{table}[]
\centering
\caption{Summary of the studies with meta-heuristic approach for the problem}
\label{tab:heuristic}
\begin{tabular}{lll}
\hline
Paper               & Solution                             & Year \\
\hline
\citet{tate1995unequal}          & GA and Penalty function  & 1995  \\
\citet{tavakkoli1998facilities}     & GA                  & 1998 \\
\citet{li2000genetic}    & GA          & 2000  \\
\citet{logendran2006methodology}    & Tabu search Algorithm     & 2006  \\
\citet{liu2007sequence} & Finding initial Solution by GA for MIP & 2007  \\
\citet{scholz2009stats}     & Tabu Search with Slicing tree       & 2009 \\   
\citet{komarudin2010applying}     & Ant Colony with Tabu Search         & 2010 \\     
\citet{scholz2010extensions}     & Tabu Search with Slicing tree         & 2010 \\                      
\hline
\end{tabular}
\end{table}

\citet{jithavech2010simulation} proposes a simulation method to solve the equal area FLP in stochastic conditions and uses a Quadratic Assignment Problem (QAP) formulation. They consider stochastic distribution for flow between each two departments and find a layout with minimum risk for a single period problem. Their algorithm minimizes the maximum expected increase in material handling cost of the layout according to their definition of risk. 

Responding to the aforementioned potentials for the short-lived electronic and food industry, \citet{rosenblatt1986dynamics} introduces a new type of FLPs. They assume predetermined periods in which the change in the layout and flow between the departments is required. The rearrangement cost is also a given value for each period. The strategic view towards the FLPs introduced the dynamic facility layout problem (DFLP) in which the material flow volumes between different departments change in certain time periods \citep{mckendall2010heuristics}. Different studies have been conducted after \citet{rosenblatt1986dynamics} that present a dynamic programming for solving equal area DFLPs (for a review see \citep{kulturel2007approaches} and \citep{drira2007facility}). 
In all of these studies, number of periods and flows in each period is known and algorithms try to seek quality solutions in all of the periods. Some constraints such as budget and unequal area departments have been added in some of the works. Some important studies of the DFLP in such categories are shown in Table \ref{tab:dynamic}.

\begin{table}[]
\centering
\caption{Summary of the multi-period dynamic flows studies approach for the problem}
\label{tab:dynamic}
\begin{tabular}{lll}
\hline
Paper               & Solution                             & Year \\
\hline
\citet{baykasouglu2001simulated}          & Simulated Annealing  & 2001  \\
\citet{balakrishnan2003hybrid}     & GA                   & 2003 \\
\citet{dunker2005combining}    & GA and Dynamic Programming            & 2005  \\
\citet{baykasoglu2006ant}    & Ant Colony with Budget constraint     & 2006  \\
\citet{mckendall2006hybrid} & Hybrid Ant System & 2006  \\
\citet{balakrishnan2006note}     & GA and Simulated Annealing       & 2006 \\   
\citet{kulturel2007approaches}     & Review dynamic and stochastic FLP         & 2007 \\     
\citet{drira2007facility}     & A survey in Facility Layout problem         & 2007 \\    
\citet{balakrishnan2009dynamic}     & Forecast Uncertainty and Rolling Horizons       & 2009 \\    
\citet{mckendall2010heuristics}     & Tabu Search and Boundary Search       & 2010 \\    
\citet{csahin2010simulated}     & Simulated Annealing with Budget limits         & 2010 \\                      
\citet{pillai2011design}     & Robust Optimization Approach         & 2011 \\  
\citet{xu2012multi}     & Fuzzy Environment with Swarm Optimization         & 2012 \\  
\citet{chen2013new}     & A new data set and Ant colony Optimization         & 2013 \\ 
\citet{mazinani2013dynamic}	&	GA	& 	2013 	\\ 
\citet{pourvaziri2014hybrid} & GA& 2014 \\
\citet{ulutas2015dynamic} & Clonal Selection Based Algorithm	& 2015 \\ 
\hline
\end{tabular}
\end{table}

However, with the need for rapid changes in facilities in order to compete and survive in the market, firms need to introduce new products to, and remove some of the old products from the production lines, which often results in changes in old facilities and usually some revisions in whole layout. To respond to this situations the new generation of DFLPs emerged, which considers changes in the type and the number of facilities as well as the material flow changes between them. Similar to the classic approach of DFLP, change in the location of each facility generates a cost. Solution to such problems is in the form of several layouts proportional to the number of time periods. As mentioned, \citet{dong2009shortest} introduce the first study in this area by considering flow between machines and rearrangement costs as given values. Representing the problem as a graph and converting it into a shortest path problem, they use a simulated annealing algorithm to solve their problem. In spite of such progress towards the DFLPs, to the best of our knowledge, we are the first to consider an unequal area DFLP with stochastic material flow, cost, and time of change. 

\section{Model Formulation}
\label{sec: model}

The result of DFLP is a layout, representing the best arrangement of the corresponding departments in the layout. To represent this arrangement, we use the slicing tree method introduced by \citet{otten1982automatic}, and frequently used in the literature \citep{diego2009solving, scholz2009stats, scholz2010extensions, wong2010comparison}. Since the FLP has proven to be an NP-hard problem, it is impossible to solve it for real world large scale problems in polynomial time. We use GA to solve each problem in order to find a near optimal solution. Since our problem parameters including material flows and costs have stochastic nature, we use simulation to evaluate the quality of each layout. GA needs deterministic parameters; however, in a dynamic setting as our paper they are not so that to attain suitable solutions we divide each problem into several static FLPs and solve them separately. Our assumptions are as followings: 
\begin{enumerate}[(a)]
\item In an indeterminable time, new departments could be added (removed) to (from) the existing layout. The number of departments to be added (removed) is known.
\item The material flow between departments will be estimated as a distribution. We assume this distribution is known. Same setting is applied to the rearrangement cost of the departments.
\item A known life cycle is determinable for the product.
\end{enumerate}

Literature review shows that normal and uniform distributions are more common than others for the setting closest to ours \citep{tavakkoli1998facilities, jithavech2010simulation} to model the flow between departments and the rearrangement costs in the simulation model. In predicting the demand for each product through market research and polls, typically an upper and a lower bound for product demand can be obtained. Therefore, with a good estimate, often this distribution can be assumed to be uniform.

\subsection{Slicing Tree}
In the slicing tree method, each layout is represented as triple row matrices \cite{otten1982automatic}. In the first row, the department sequence, a sequence of $n$ departments is included to assign. In the second row a sequence of $n-1$ points are included for slicing, and in the third row, for each one of $n-1$ slicing points, an orientation (either vertical or horizontal) is included. The slicing sequence matrix has a sequence of numbers from $1$ to $n-1$. In each position of the third row, a number is entered that is either $0$ or $1$, where $0$ stands for a horizontal slice and $1$ for a vertical slice. Therefore, each layout can be shown as a $3\times n$ matrix, and for a problem with $n$ departments $2^{n-1}n!(n-1)!$ layouts are possible. However, the output of the slicing tree would not guarantee the condition that each department satisfies the maximum ratio of length to width; so, another algorithm should be built to satisfy this condition. 

\subsection{Genetic Algorithm}
\begin{figure}
\begin{center}
\includegraphics[width=4in]{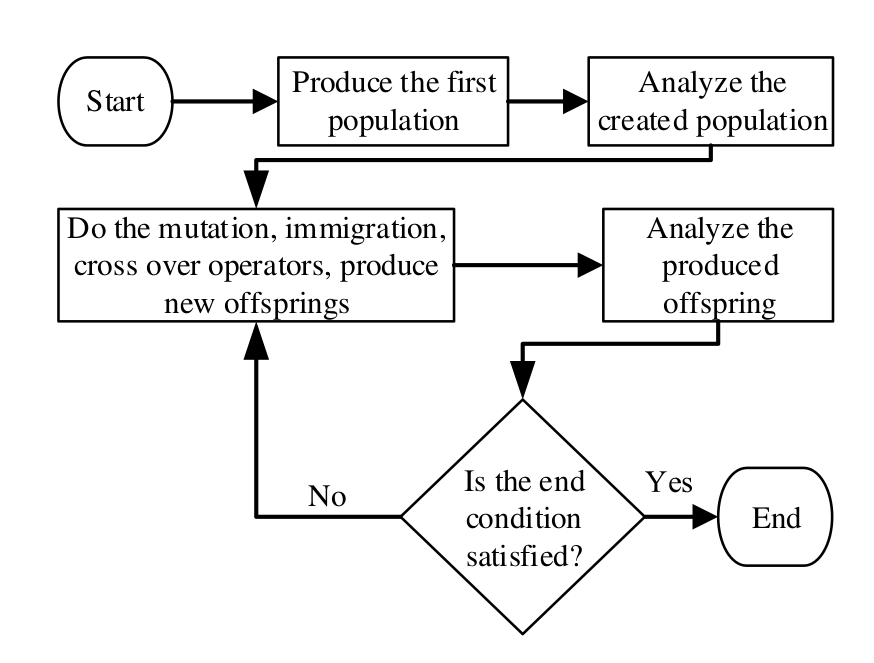}
\end{center}
\caption{The designed GA}
\label{fig:GA}
\end{figure}

GA was introduced based on inheriting characteristics by genes and is categorized as a random search. GA is used for solving various FLPs, among which \cite{azadivar2000facility}, \cite{wu2002optimisation}, and \cite{dunker2003coevolutionary} investigate deterministic FLPs and \cite{balakrishnan2003hybrid} and \cite{dunker2005combining} consider DFLPs. Our customized GA is as follows:

We use a $3\times n$ matrix to demonstrate each layout as a chromosome. Each population comprises some of these chromosomes. To produce each population, we use the best solution from the last population. Our customized GA is shown in Figure \ref{fig:GA}. 

The algorithm will continue until a given number of generations without improvement in objective function is reached. Details of each step are explained in following sections.

\subsubsection{Initial Solution}
We observed that finding a promising initial solution has a very positive effect on the performance of genetic algorithm. Our GA starts from within the feasible area and improves the solutions iteratively. Consequently, the probability of finding better or near optimal solutions increases. The best initial solution could be the one in which the departments with larger flows in between compared to other departments are placed nearby \citep{meller1996facility}. In slicing tree representation, placing the departments near each other in department sequence, creates no necessity for them to be nearby in the final sequence. Therefore to place two departments with big flows nearby, it is necessary for them to be placed at the end of the sequence in slicing time. Considering this idea, we propose the following two methods for putting two departments next to each other. 

\begin{enumerate}[(1)]
\item Put the two corresponding departments at the end of the department sequence and put the number for slicing them at the end of the slicing sequence. Under this condition the two departments will be nearby in the final layout. 

\item Put the two departments at the beginning of the department sequence. Put numbers 1 and 2 in the slicing sequence in the order of 2 and then 1. Following this rule for their order, 2 and then 1, you can put them nearby at any position in the slicing sequence. This condition will guarantee the fact that they will be nearby in the ultimate layout. This is due to the fact that in the previous slices, regardless of their sequence, no slice has been done between the two departments. The first slice occurs when one gets to number 2. This slice will separate the two departments as a cubic rectangular mass. In the next step the slice with number 1 will divide this mass into two parts, one for each department. Therefore, the two departments will be nearby in the final layout. 
\end{enumerate}

We use these two methods and propose the following algorithm to populate the initial solutions consisting five steps as follows:

\begin{enumerate}[(a)]
\item Select two departments with maximum amount of flow. Put them at the end of department sequence. Put number $n-1$ at the end of slicing sequence. 
\item Select the next pair of departments with maximum amount of flow and place them at the start of the department sequence. 
\item Select position $\hat{n}$ from 1 to $n-2$ randomly. If $\hat{n}\neq n-2$ then go to step (d), otherwise go to step (e). 
\item Put number 2 in position $\hat{n}$ and number 1 in position $\hat{n}+1$ in the slicing sequence resulted from step (b).
\item Create a binary number randomly. If the number created is 1, then put number 2 in position $n-3$ and number 1 in position $n-2$ in the slicing sequence resulted from step (b), otherwise, put number 2 in position 1 and number 1 in position 2 of the slicing sequence. 
\end{enumerate}

Running this algorithm can potentially reach to $2^{n-1}(n-4)!(n-4)!(n-1)$ for possible solutions which is significantly lower than the previous value of $2^{n-1}n!(n-1)!$. Table \ref{tab:ST} compares the number of possible solutions with and without the initial solution.
As shown, number of solutions is considerably reduced in our algorithm, which in middle size problems results in an invaluable computation time reduction.

\begin{table}[]
\centering
\caption{Number of possible solutions for a problem with slicing tree representation method}
\label{tab:ST}
\begin{tabular}{ccc}
\hline
Number of               & Number of       & Number of possible solutions  \\
departments               & possible solutions      & after creating an initial solution \\
\hline
7          & 232243200 & 13824  \\
8    & 26011238400  & 516096 \\
9   & 3.74562E+12 & 29491200  \\
10    & 6.74211E+14  & 2388787200  \\
11 & 1.48326E+17 & 2.60112E+11  \\
12     & 3.91582E+19 & 3.66238E+13 \\   
20    & 1.55163E+41 & 4.36077E+33 \\     
30     & 1.25912E+72 & 2.53225E+63 \\    
35     & 5.24104E+88 & 3.9495E+79 \\    
62     & 3.6832E+187 & 7.7715E+176 \\    
\hline
\end{tabular}
\end{table}

\subsubsection{Crossover}
In each iteration of GA, a new population is produced. One of the operators used to create a new population is the crossover operator. This operator selects two chromosomes and with transmission of the characteristics from both, two offspring will be produced, which have some of their characteristics from the first parent and some of them from the second parent.

One point and two point cross overs are the most well-known cross overs and used here. In one point cross over, a point is selected randomly and the genes before the point will be directly transmitted to the produced offspring. The genes after the selected point are transmitted from parent 1 to offspring 2 and vice versa. To make the offspring acceptable, changes are necessary as follows: in the correction process we define the unmentioned genes in each offspring and assign them based on their sequence in the related parent right after the crossing gene. Also, we remove the duplicated genes that are placed after the selected gene. At the end, we put the other genes with a sequence same as parent's sequence at the end of offspring. (see Tables \ref{tab:OC1}, \ref{tab:OC2}, and \ref{tab:OC3})

\begin{table}[H]
\centering
\caption{One point Cross over selecting parents}
\label{tab:OC1}
\begin{tabular}{|c|c|c|c|c|c|c|c|c|c|}
\hline
Parent 1: & {\cellcolor[gray]{.8}}8 & {\cellcolor[gray]{.8}}4 & {\cellcolor[gray]{.8}}2 & {\cellcolor[gray]{.8}}6 & 7 & 3 & 9 & 1 & 5\\
\hline
Parent 2: & {\cellcolor[gray]{.8}}2 & {\cellcolor[gray]{.8}}9 & {\cellcolor[gray]{.8}}5 & {\cellcolor[gray]{.8}}8 & 3 & 4 & 6 & 7 & 1\\
\hline
\end{tabular}
\end{table}

\begin{table}[H]
\centering
\caption{One point Cross over producing offspring}
\label{tab:OC2}
\begin{tabular}{|c|c|c|c|c|c|c|c|c|c|}
\hline
Offspring 1: & \cellcolor[gray]{0.8}8 & \cellcolor[gray]{0.8}4 & \cellcolor[gray]{0.8}2 & \cellcolor[gray]{0.8}6 & 3 & 4 & 6 & 7 & 1\\
\hline
Offspring 2: & \cellcolor[gray]{0.8}2 & \cellcolor[gray]{0.8}9 & \cellcolor[gray]{0.8}5 & \cellcolor[gray]{0.8}8 & 7 & 3 & 9 & 1 & 5\\
\hline
\end{tabular}
\end{table}

\begin{table}[H]
\centering
\caption{One point Cross over correcting the sequence of produced offspring}
\label{tab:OC3}
\begin{tabular}{|c|c|c|c|c|c|c|c|c|c|}
\hline
Child 1: & \cellcolor[gray]{0.8}8 & \cellcolor[gray]{0.8}4 & \cellcolor[gray]{0.8}2 & \cellcolor[gray]{0.8}6 & 9 & 5 & 3 & 7 & 1\\
\hline
Child 2: & \cellcolor[gray]{0.8}2 & \cellcolor[gray]{0.8}9 & \cellcolor[gray]{0.8}5 & \cellcolor[gray]{0.8}8 & 4 & 6 & 7 & 3 & 1\\
\hline
\end{tabular}
\end{table}

In two point cross over, two points are selected randomly. The genes that are between these points are transmitted to all offspring, and the rest will be transmitted from parent 1 to offspring 1 and from parent 2 to offspring 2. To make the offspring acceptable, the required revisions are made and the new population will be produced. The process is shown in Tables \ref{tab:TC1}, \ref{tab:TC2}, and \ref{tab:TC3}.

\begin{table}[H]
\centering
\caption{Two point Cross over selecting parents}
\label{tab:TC1}
\begin{tabular}{|c|c|c|c|c|c|c|c|c|c|}
\hline
Parent 1: & 1 & 9 & 7 & 5 & {\cellcolor[gray]{.8}}6 & {\cellcolor[gray]{.8}}4 & {\cellcolor[gray]{.8}}2 & 3 & 8\\
\hline
Parent 2: & 6 & 5 & 3 & 9 & {\cellcolor[gray]{.8}}8 & {\cellcolor[gray]{.8}}4 & {\cellcolor[gray]{.8}}7 & 2 & 1\\
\hline
\end{tabular}
\end{table}

\begin{table}[H]
\centering
\caption{Two point Cross over producing offspring}
\label{tab:TC2}
\begin{tabular}{|c|c|c|c|c|c|c|c|c|c|}
\hline
Offspring 1: & 1 & 9 & 7 & 5 & {\cellcolor[gray]{.8}}8 & {\cellcolor[gray]{.8}}4 & {\cellcolor[gray]{.8}}7 & 3 & 8\\
\hline
Offspring 2: & 6 & 5 & 3 & 9 & {\cellcolor[gray]{.8}}6 & {\cellcolor[gray]{.8}}4 & {\cellcolor[gray]{.8}}2 & 2 & 1\\
\hline
\end{tabular}
\end{table}

\begin{table}[H]
\centering
\caption{Two point Cross over correcting the sequence of produced offspring}
\label{tab:TC3}
\begin{tabular}{|c|c|c|c|c|c|c|c|c|c|}
\hline
Child 1: & 1 & 9 & 5 & 6 & {\cellcolor[gray]{.8}}8 & {\cellcolor[gray]{.8}}4 & {\cellcolor[gray]{.8}}7 & 2 & 3\\
\hline
Child 2: & 5 & 3 & 9 & 8 & {\cellcolor[gray]{.8}}6 & {\cellcolor[gray]{.8}}4 & {\cellcolor[gray]{.8}}2 & 7 & 1\\
\hline
\end{tabular}
\end{table}

\subsubsection{Mutation}
The Mutation operator works on one parent solution and generates an offspring by modifying the parent solution's features in random conditions. It is useful to create diverse population and also escape from local optima. In our mutation operator, we choose two random numbers $N_1$ and $N_2$ and exchange their position, where $2<N_1<N_2<n-1$. Table \ref{tab:MO} represents an instance, which considers a sample chromosome and numbers 2 and 6 are randomly chosen to be exchanged. 

\begin{table}[H]
\centering
\caption{Mutation operator}
\label{tab:MO}
\begin{tabular}{|c|c|c|c|c|c|c|c|c|c|}
\hline
Parent solution before mutation: & 4 & 8 & {\cellcolor[gray]{.8}}6 & 3 & 9 & {\cellcolor[gray]{.8}}2 & 7 & 5 & 1\\
\hline
Offspring made by mutation operator: & 4 & 8 & {\cellcolor[gray]{.8}}2 & 3 & 9 & {\cellcolor[gray]{.8}}6 & 7 & 5 & 1\\
\hline
\end{tabular}
\end{table}

\subsubsection{Migration}
\begin{figure}
\begin{center}
\includegraphics[width=2in]{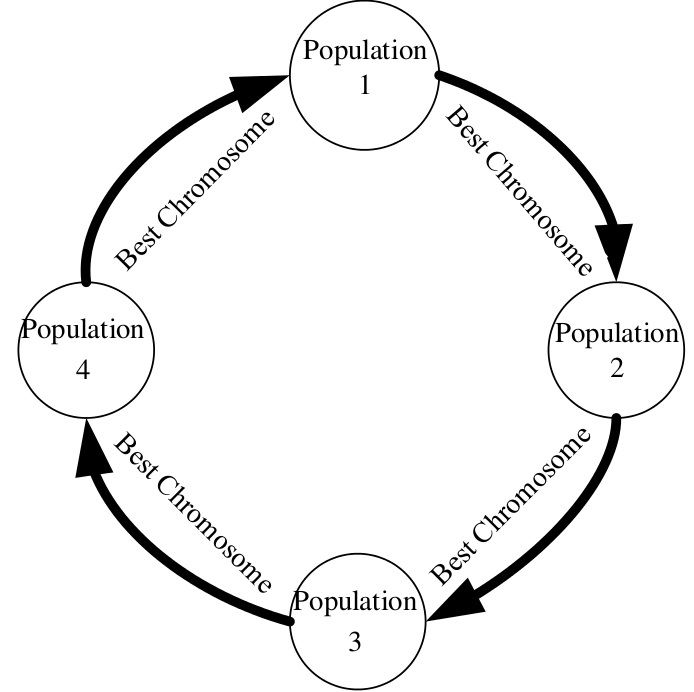}
\end{center}
\caption{The best instances in a population migrate to the next neighbour in the ring}
\label{fig:Mig}
\end{figure}

In the process of the proposed GA, four populations are created and the GA runs in each population separately. At the end of each run, the best instances of each population migrates to the next population. This procedure is performed iteratively to insure steady improvement. Figure \ref{fig:Mig} represents this procedure.

\subsubsection{Dynamic calculation of the number of operators }
In our algorithm in each stage of the proposed GA, the number of offspring that are produced from each operator is determined in a unique way. In each iteration, considering the percentage of improvement of the solution in the previous iterations, the number of operators is selected so that it could help getting quicker improvement. \citet{kianfar2012study} used the last improved percentage of the objective function to determine the portion of each operator in the next iteration. We suggest using an average of improvement percentages in the five past iterations which makes a more reliable and stable basis to get better results and escape local optima. In those iterations that the improvement percentage is zero, it is probable that the algorithm is in local optima. In order to escape from this condition, it is required to increase the number of random operators to allow more explorations. Thus, we increase the number of mutation and migration operators in our GA. Also, in iterations that the improvement percentage is suitable enough, the number of cross over operators increases and the random operations decrease.

\subsubsection{Tuning HGA parameters} 
Parameter tuning is one of the most important steps in building a meta-heuristic algorithm. We have five parameters as {\em initial population size, number of steps without improvement for stopping criteria, dynamic child parameters of crossover, and mutation and migration}. The tuning process that was introduced by \citet{kianfar2012study} is appropriate for our setting too. The process is that a selected number of parameters are fixed in certain amounts and other parameters change within a given boundary. Considering the value of objective function, the best value for each parameter is selected. 

With experimental tests to choose the best value for parameters {\em initial population size} and {\em number of steps without improvement}, the amount of 70 and 300 are selected respectively. Table \ref{tab:Tuning} demonstrates the tuned amounts of {\em dynamic child parameters of crossover, mutation parameter and migration parameter} based on the best average of five previous improvement percentages in objective function. 
In this table, $Impr$ represents the average percentage of the objective function in five previous iterations and the numbers included in each column are the percentage of offspring produced by each operator. In each stage the amount of $Impr$ is calculated first and using it, the number of operators is calculated and the process continues.

\begin{table}[]
\centering
\caption{Percentage of produced offspring based on the percentage of growing in answer in last stage}
\label{tab:Tuning}
\begin{tabular}{cccc}
\hline
Improvement  & Crossover       & Mutation & Migration  \\
\hline
0&				0.61&	0.31&	0.08 \\
$0<Impr<1$&		0.67&	0.27&	0.06\\
$1\leq Impr<2$&	0.77&	0.19&	0.04\\
$2\leq Impr<4$&	0.80&	0.15&	0.05\\
$4\leq Impr<6$&	0.87&	0.10&	0.03\\
$6\leq Impr<8$&	0.89&	0.08&	0.03\\
$8\leq Impr$&	0.92&	0.05&	0.02\\
\hline
\end{tabular}
\end{table}

\subsubsection{Objective Function}
To evaluate the efficiency of each solution we use the total amount of material handling cost for objective function. On the other hand, the solutions created by slicing tree are not guaranteed to be feasible. To deal with this issue the  we use the penalty objective function as proposed by \citet{komarudin2010applying}. In this objective function, if a department violates the limitation of the maximum ratio of the length to width, an amount is added as the penalty to the objective function. Using this objective function, GA tries to produce solutions that can guarantee the limitations. Our objective function is as follows: 

\[\sum\limits_{j=1}^n \sum\limits_{i=1, i\neq j}^nf_{ij}d_{ij} + p_{inf}(V_{feas} - V_{all}) + \sum\limits_{i=1}^n Re_i \times ReCost_i,\]
in which the first part stands for the total amount of material handling cost. $p_{inf}$ is the number of infeasible departments. $V_{feas}$ is the best amount calculated for a layout in which all the departments are feasible, $V_{all}$ is the best amount calculated for a layout regardless of its feasibility, $Re_i$ is 1 if department $i$ is rearranged and 0 if not, and $ReCost_i$ is the rearranging cost of department $i$.

%The objective function of the model is not a dynamic equation like the objective functions in \cite{rosenblatt1986dynamics}, \cite{balakrishnan2003hybrid}, \cite{kulturel2007approaches}, and \cite{drira2007facility}. The dynamicity of our model is conducted in the material flow parameters. % Furthermore, in our problem the rearrangement cost is not a function of distance as \citet{mckendall2010heuristics}.

%\section{Solution Approach} \label{sec: solution}

\subsection{Hybrid simulation for DFLP} \label{sec: solution}

Emerging new products usually results in adding new machines in the layout. Also, because the design of layout should be conducted with respect to the marketing data for the new product, the material flows are not exact as the marketing data are not exact. However, we can define the material flow within lower and upper bonds as a uniform distribution. 

On the other hand, our search engine, i.e. GA needs deterministic parameters and it is not possible to use it with random parameters. The common approach to deal with this kind of problems is considering material flows as $a\mu + b\sigma$ with a fixed $a$ and $b$, e.g. $\mu + \sigma$ and using it in the deterministic algorithms (see \citep{tavakkoli2007design, jithavech2010simulation} ).
However, this is not the best possible approach and there is no guarantee on the solution's quality in the real situation. In other words, utilizing another deterministic material flow such as $\mu + 2\sigma$ may result in better solution than that of $\mu + \sigma$. Therefore, searching among larger sets of $a$ and $b$ may result in better solutions.

On the other hand, in continuous environments there are countless points and we cannot evaluate all possible sample material flows. To address this issue, we consider a few number of initial material flows and with the help of an iterative greedy search algorithm we find a better solution. In the greedy search, for each selected material flow we calculate the corresponding cost and then try to move to a direction which decreases the cost. In this process, one must determine a criterion to compare the resulted layouts. Note that, the classic objective function, $f_{ij} d_{ij}$, is not an appropriate measure here. For example, consider obtained equal layouts $X$ and $Y$ corresponding to material flows of $\mu + \sigma$ and $\mu - \sigma (>0)$. Considering the material flows, the corresponding objective function of layout $X$, surely is greater than that of $Y$. This is due to the fact that the two layouts are the same and as a result they have equal $d_{ij}, ~ \forall i \neq j \in \{1,\dots, n\}$. Considering the objective function, i.e., $\sum \limits_{i \neq j}^n f_{ij} d_{ij}$, problem $X$ has a greater objective function than $Y$. Therefore, layout $Y$ will be the choice among these two, which may not be the right decision.

Therefore, in order to compare different layouts, a measuring framework is required, not sensitive to the corresponding material flow which is used in objective function of design procedure. 
Such framework should compare all layouts with a same material flow to obtain comparable costs. Simulation ensures such criteria and can be used to measure how good the obtained layouts are. 
Thus, we provide simulation which gets a certain number of layouts---with different amounts of material flow (different $a$ and $b$ in $a\mu + b\sigma$)---and analyze them with a same material flow to provide the corresponding objective function of each layout.

\begin{figure}
\begin{center}
\includegraphics[width=6in]{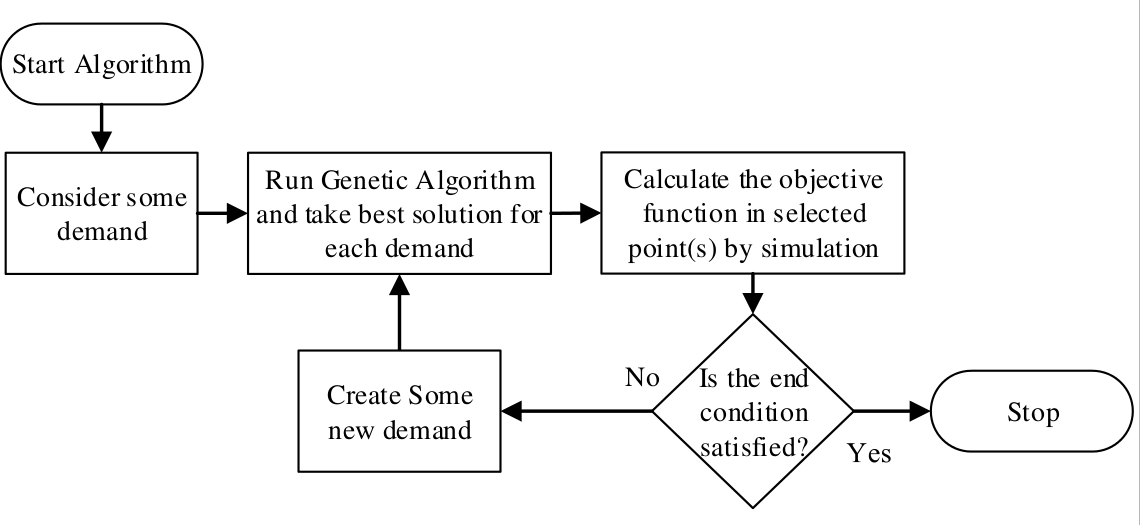}
\end{center}
\caption{The algorithm used to optimize the solution}
\label{fig:OptAlg}
\end{figure}

Our algorithm to solve the DFLP with stochastic material flow, stochastic cost, and stochastic time of change is presented in Figure \ref{fig:OptAlg}. 
The algorithm gets $k$ material flows (set of $(a, b)$) and for each the best possible layout is obtained through GA. Then, objective value of all layouts are calculated by the simulator, which its details is provided in section \ref{sec:simulator_details}. %parameters of which are set by these sets of $a$ and $b$. 
An ANOVA experiment analyzes the resulted cost (objective function values), to sort the layouts according to the their costs---if they are statistically different. 
Then, the worst layout and the corresponding material flow will be removed and a new material flow according to the algorithm will be generated. This procedure continues until the stopping condition is met. 

In each cycle of simulation, all possible pairs of layouts are compared with each other. Among statistical comparison methods 
%LSD, Scheffe, Duncan, Tukey, and Orthogonal Contrast methods \citep{montgomery2008design}, we chose Tukey test, which analyzes the equality hypothesis for the cost of layouts. Rejecting the null hypothesis results in a new flow and algorithm will continue.
we chose Tukey test \citep{montgomery2008design}, which analyzes the equality hypothesis for the cost of layouts. Rejecting the null hypothesis results in a new flow and algorithm will continue.
Tukey test rejects the null hypothesis of equality of costs of two material flows if $c_i-c_j>\frac{q_a(p,f)}{\sqrt{2}}\sqrt{MSE(\frac{1}{n_i}+\frac{1}{n_j})}$, where $c_i$ represents the cost of $i^{th}$ material flow, $q_a(p,f)$ is the Tukey test factor, MSE is mean square error and $n_i$ is the number of observations in $i^{th}$ material flow. 
After getting the results of the pairwise comparison, if the null hypothesis is rejected the material flow that has the worst cost compared to other material flows will be removed. %This means that the material flow that Tukey test rejects its null hypothesis will be removed. 
Regarding the minimization nature of the objective function, the corresponding layout of this material flow has the highest objective function value among all layouts. In order to substitute a new material flow with the removed material flow, we select the two corresponding material flows of the layouts that have sufficiently good results when compared with other material flows. The average of these material flows will be used in the next cycle. This algorithm will continue until one of the following conditions is satisfied:

\begin{enumerate}[(1)]
\item Simulation takes longer than a given time. 
\item Equality of layouts corresponding to different material flows is not rejected by Tukey test.
\end{enumerate}

\subsubsection{Simulator}\label{sec:simulator_details}
\begin{figure}
\begin{center}
\includegraphics[width=6in]{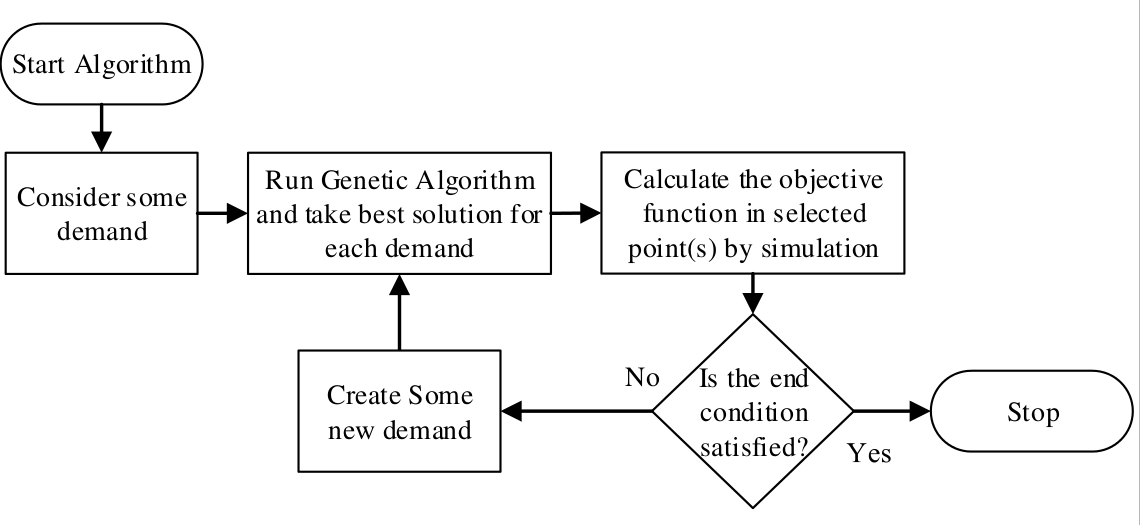}
\end{center}
\caption{Each run of simulation}
\label{fig:Sim}
\end{figure}

In our simulation algorithm, first a random uniform demand is created, defining the flow between each pair of departments. Then, the distance between the departments pairs is calculated and considering the material flow between them, the material handling cost is calculated. Then algorithm goes back to the initial step and next random demand will be created. This iterative algorithm repeats until the stopping criteria is met. In each run for simulation the process in Figure \ref{fig:Sim} is used. To ensure the accuracy of the results, the simulation has been repeated 10000 times for each layout plan.

\section{Computational Experiments}
\label{sec: computational}

In this section the results of numerical experiments is presented to show the effectiveness of the GA and its initial population. Then, the result of the hybrid algorithm to deal with DFLP is presented. 

\subsection{Performance of The GA algorithm}
%GA, the solution engine of our algorithm is very important and must be reliable. As our GA is linked with the corresponding simulation algorithm which might run for extensive amount of rounds to ensure accuracy, the computation time of GA becomes important. Therefore, shorter execution times for GA while ensuring high quality solutions is of great importance. 
In order to examine the performance of our suggested GA, we evaluate it with some benchmark problems, selected from \cite{komarudin2010applying}. The maximum aspect ratio of each department for problems 07, 08, 09 and AB20 is 4 and for SC30 is 5. 
The results of our GA algorithm are summarized in Table \ref{tab:Comp}, which as it is shown our GA provides the optimal solution in small problems. 
In our algorithm, the GA continues running until one of the following conditions is satisfied:

\begin{itemize}
\item The number of times in which a better solution has not been observes surpasses 300.
\item The total number of runs for the algorithm surpasses 1000.
\end{itemize}

\begin{table}[H]
\centering
\caption{Comparison between the output of proposed algorithm and and studies existing in litrature}
\label{tab:Comp}
\begin{tabular}{cccc}
\hline
Problem  & Best In \citet{komarudin2010applying} & Our Solution & Gap  \\
\hline
07&	131.68&	131.68&	0\\
08&	243.16&	243.12&	0\\
09&	239.07&	239.05&	0\\
AB20&	5225.96&	5543.5&	0.09\\
SC30&	3707&	4248&	0.146\\
\hline
\end{tabular}
\end{table}
%
%To evaluate the reliability of the algorithm for each example, the algorithm is repeated five times. The complete results including the average run time, the number of populations produced for running the algorithm, and the initial population for each example are shown in Table \ref{tab:Comp2}.
%%
%\begin{table}[H]
%\centering
%\caption{Comparison between the output of proposed algorithm and and studies existing in litrature}
%\label{tab:Comp2}
%\begin{tabular}{cccccccccc}
%\hline
%Problem&	Best &	Worst &	Average  &	Solution &	\citet{komarudin2010applying} &	Time &	Population &	Turn\\
%&	Result&	 Result&	 Result &	 Time (Avg)&	 Best Time&	 Improvement&	 Number&	\\
%\hline
%07&	131.68&	137.05&	132.754&	{\cellcolor[gray]{.7}}65&	{\cellcolor[gray]{0.9}}756&	91\%&	120&	175\\
%08&	243.11&	254&	247.466&	{\cellcolor[gray]{.7}}87&	{\cellcolor[gray]{0.9}}792&	89\%&	120	&250\\
%09&	239.05&	254&	242.04&	{\cellcolor[gray]{.7}}114&	{\cellcolor[gray]{0.9}}900&	87\%&	120	&250\\
%AB20&	5543.5&	5847&	5700.2&	{\cellcolor[gray]{.7}}3795&	{\cellcolor[gray]{0.9}}17820&	79\%&	125	&1000\\
%SC30&	4248&	4768&	4609.4&	{\cellcolor[gray]{.7}}6129&	{\cellcolor[gray]{0.9}}23544&	74\%&	125	&1000\\
%\hline
%\end{tabular}
%\end{table}
%
%It is notable that the average time for each problem is considerably reduced compared to the \citet{komarudin2010applying} model, which is the direct effect of good initial solutions. The Average time is reduced at least by 74\%. The resulted layouts are presented in Appendix A.

\subsection{Performance of Hybrid GA}
Since our problem is unprecedented in the literature, there was no study and publicly available dataset in dynamic conditions for benchmarking. To evaluate the performance of the algorithm, benchmark problem 09 which has 9 departments is used. We considered the case that one department is added to the layout and the other departments remained fixed, which is a common case in food or electronic industries. 
Also, since we did not have given value for rearrangement costs and the maximum ratio for each department, we check different values for them to analyze their effects.
This rearrangement cost affects the final layout in a way that by increasing this cost, GA will converge to a layout without any change and non-rearranged departments and vice versa.
Since there is no given value from real world cases, we considered an upper and lower bound for rearrangement cost and for three different intervals (0, 0), (110, 130) and (130, 150) examined the solutions. For each problem we generated a random rearrangement cost using the uniform distribution. Distribution (0, 0) leads to a layout without considering the rearrangement cost. Distribution (130, 150) is related to the material flow in which all experiments lead to a layout without any changes. With considering material handling cost of the layout in which the new department is added without any changes, this uniform distribution (130, 150) was obtained. Lastly, (110, 130) is the cost that can moderate the rearrangement cost and the material handling cost. It is noticeable that typically in real world the rearrangement cost is deterministic.

Similarly, we checked two values for the maximum ratio of each department. By definition of problem 09, maximum ratio for each department is 4, and for the new department also the ratio is set to 4. This ratio restricts the layout and obligates it to rearrange some of the departments, and algorithm cannot add the new department in the new space of the layout without any changes (the maximum ratio penalty does not allow it). Thus, we also considered maximum ratio 13 for the new department to see the effect of having large maximum ratio. This amount comes from the ratio of height to width of improvement area which is added to the layout. In other words, the ratio allows new department to be added to empty space without any change in the initial layout.

%As mentioned before, the problem is solved for several set of material flows and, considering the best results, the next amount for the material flow is selected and this process continues to obtain a better layout. For each amount, material flow between the new department and each one in the existing layout is produced randomly by a uniform distribution. The upper and lower bounds of the uniform distribution are used as the maximum and minimum amount of material flow for this department respectively. 
Finally, according to Algorithm \ref{fig:OptAlg}, we need to have some initial material flows which we consider:
$$[\mu - \sigma,\mu, \mu+\sigma, \mu+1.5\sigma, \mu+2\sigma],$$
in which $\mu + 1.5\sigma$ is considered as \citet{tavakkoli2007design} suggests. For each amount of material flow, the best layout of proposed GA is used and then the simulation provides the cost of the layout. %Then, the new material flow is selected considering the costs corresponding to each material flow and it is fed into the algorithm. 

\begin{figure}
\centering
\includegraphics[width=5in]{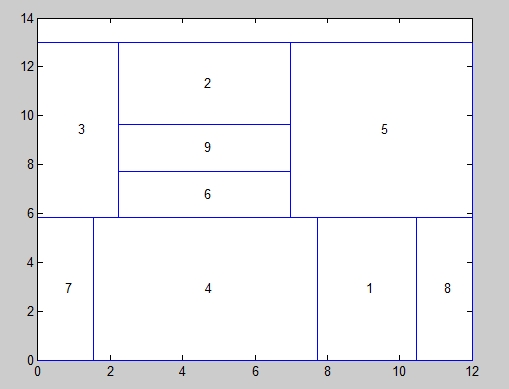}
\caption{Initial layout for the dynamic problem}
\label{fig:Lay05}
\end{figure}

The combination of these three parameters results in 30 different problems, e.g. the problem that the material flow is $\mu$, rearrangement cost is a uniform distribution of (110, 130) and the maximum ratio of the new department is 4. %In each problem five tests are performed.
We examine the effect of flow, rearrangement cost, and maximum ratio, and then find the best layout among all. The initial layout for 9 departments example is shown in Figure \ref{fig:Lay05} and the subsequent figures and tables show the effect of different parameters. 

In order to analyze the result, we designed a factorial experiment and the result of ANOVA is shown in Table \ref{tab:ANOVA}. As shown, two parameters, material flows and rearrangement cost are statistically different and null hypothesis of equality of their costs and the corresponding material flows are rejected. Also, the null hypothesis of equality for the effect of maximum ratios on the solution cannot be rejected. 
Thus, for each rearrangement cost we run the hybrid GA for different material flows and the corresponding objective values are reported in sections \ref{sec:result_0_0}, \ref{sec:result_110_130}, and \ref{sec:result_130_150} for rearrangement costs (0,0), (110,130), and (130,150) respectively.

\begin{table}
\centering
\caption{ANOVA results for the experiment in dynamic condition}
\label{tab:ANOVA}
\begin{tabular}{cccccc}
\hline
Sum of Squares&     &		df&	Mean Squares&	F-ratio&	P-Value \\
\hline
SSflow&	4.23117E+11&	4&	1.05779E+11&	23.76275&	0.0001\\
SScost&	9.45034E+11&	2&	4.72517E+11&	106.1484&	0.0001\\
SSratio&	1868271248&	1&	1868271248&	0.419697&	0.5183\\
SSE&	5.34177E+11&	120&	4451472930	  & & \\	 
\hline
\end{tabular}
\end{table}

\subsubsection{Result of cost U(0 , 0)}\label{sec:result_0_0}

\begin{figure}
\begin{center}
\includegraphics[width=5in]{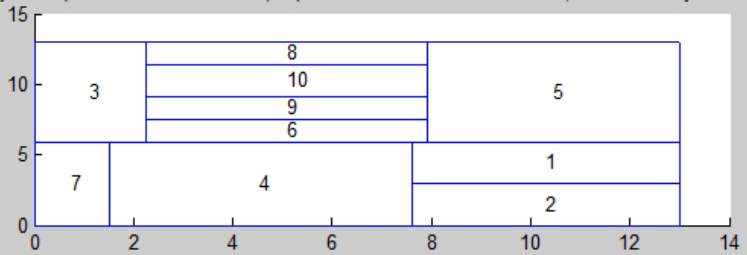}
\end{center}
\caption{Resulted test in dynamic condition with U(0,0) for rearrangement cost for all departments}
\label{fig:Lay06}
\end{figure}

The first run by five material flows is done and \emph{P-value} for the experiment is 0.31. Therefore, we cannot reject the hypothesis of equality of the material flows. Table \ref{tab:OutSim} provides the cost of each layout and Figure \ref{fig:Lay06} shows the best layout according to the simulation objective function related to $\mu - \sigma$.

In this layout, no cost is charged for changes in locations of the departments and consequently, as shown, most of the departments have changed their location compared to the previous layout (except for departments 3, 4, 7). In simulation, the cost of each layout is obtained considering the random material flows in the product life cycle as a reliable criterion for selecting the best layout. 

%$U(0 , 0)$ as the rearrangement cost has the best result of layout with 10 departments. 
\begin{table}[H]
\centering
\caption{Output of simulation for five layout plans with U(0,0) distribution for costs}
\label{tab:OutSim}
\begin{tabular}{ccc}
\hline
	& $\sum_{ij} f_{ij} d_{ij}$& 	Simulation \\
\hline
$\mu$ & 331&	3293267\\
$\mu+\sigma$ &421&	3294711\\
$\mu-\sigma$ &227&	3132879\\
$\mu+2\sigma$ &510&	3294711\\
$\mu+1.5\sigma$ &466&	3298241\\
\hline
\end{tabular}
\end{table}

\subsubsection{Result of Cost U(110, 130)}\label{sec:result_110_130}
Similarly, the algorithm is run for five initial material flows and we obtained a \emph{P-value} lower than 0.0001. According to the algorithm, pairwise comparisons for material flows were obtained and next material flow was defined. The result of pairwise comparison is shown in Table \ref{tab:PairComp}. The reference value according to Tukey formula in 95\% of significance is 44948.75. Therefore, except for three last comparisons, other material flows are different. 

\begin{table}[H]
\centering
\caption{Result of pairwise comparison by Tukey test for dynamic condition}
\label{tab:PairComp}
\begin{tabular}{ccccc}
\hline
$O_{\mu_i}$ & $O_{\mu_i}-O_{\mu_j}$	& Difference & $O_{\mu_i}-O_{\mu_j}$& Difference  \\
\hline
$O_{\mu}=3370616$ & $O_{\mu-\sigma} - O_{\mu+\sigma}$ &	 -145359.67 & $O_{\mu} - O_{\mu+\sigma}$& 90745.43  \\
$O_{\mu+\sigma}=3461362$ & $O_{\mu-\sigma} - O_{\mu}$ 	& -54614.24 & $O_{\mu} - O_{\mu+1.5\sigma}$& -54296.08  \\
$O_{\mu-\sigma}=3316002$& $O_{\mu-\sigma} - O_{\mu+1.5\sigma}$ 	& -108910.32 & $O_{\mu+\sigma} - O_{\mu+2\sigma}$& {\cellcolor[gray]{0.9}}18355.84  \\
$O_{\mu+2\sigma}=3443006$ & $O_{\mu-\sigma} - O_{\mu+2\sigma}$ 	& -127003.82 & $O_{\mu+\sigma} - O_{\mu+1.5\sigma}$& {\cellcolor[gray]{0.9}}36449.35  \\
$O_{\mu+1.5\sigma}=3424912$ & $O_{\mu} - O_{\mu+2\sigma}$ 	& -72389.59 & $O_{\mu+2\sigma} - O_{\mu+1.5\sigma}$& {\cellcolor[gray]{0.9}}18093.50  \\
\hline
\end{tabular}
\end{table}

By comparing the material flows, we remove $(\mu+\sigma)$ and assign a new material flow, which here is the average of $(\mu-\sigma)$ and $(\mu)$ and equal to $(\mu-0.5\sigma)$. It is considerable that $O_{\mu-\sigma}$ and $O_{\mu}$ are significantly lower than other material flow's objective values. In the next two iterations, the \emph{P-value} for SSmean is also lower than 0.0001 and the null hypothesis of equality of material flows is rejected. Material flows $(\mu+2\sigma)$ and $(\mu+1.5\sigma)$ are removed, and instead $(\mu-0.75\sigma)$ and $(\mu-0.625\sigma)$ are added respectively. In the 4$^{th}$ iteration, we cannot reject the hypothesis of equality of the material flows, and since the corresponding cost of $(\mu-0.625\sigma)$ is the lowest compared to others, we select this material flow as the best one and get the final result. 

\begin{figure}
\begin{center}

\end{center}
\centering
\includegraphics[width=5in]{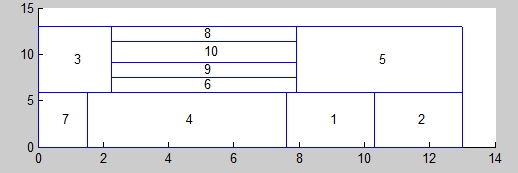}
\caption{The best resulted in dynamic condition with U(110,130) for rearrange cost}
\label{fig:Lay07}
\end{figure}

Figure \ref{fig:Lay07} is the output of the problem with a material flow of $(\mu-0.625\sigma)$. The cost of each layout is provided in Table \ref{tab:ObjRearrange}. As it is shown, the obtained results from the simulation does not lead to a same result by classic objective proposed by \cite{tavakkoli2007design}. Also, because of the high rearrangement costs, some departments have not changed their location. Departments 1, 3, 4, 6, 7, 9 are at the same locations as their previous ones.

\begin{table}[H]
\centering
\caption{The objective results for five layouts with U(110,130) for rearrange cost}
\label{tab:ObjRearrange}
\begin{tabular}{ccc}
\hline
	& $\sum_{ij} f_{ij} d_{ij}$& 	Simulation \\
\hline
$\mu$ & 1455&	3305672 \\
$\mu-0.5\sigma$ &1410&	3305672 \\
$\mu-\sigma$ &{\cellcolor[gray]{0.9}} 1354&	3206259 \\
$\mu-0.75\sigma$ &1381&	3224610 \\
$\mu-0.625\sigma$ &1382&	{\cellcolor[gray]{0.9}}3127449 \\
\hline
\end{tabular}
\end{table}

\subsubsection{Result of cost U(130, 150)}\label{sec:result_130_150}

\begin{figure}
\begin{center}
\includegraphics[width=5in]{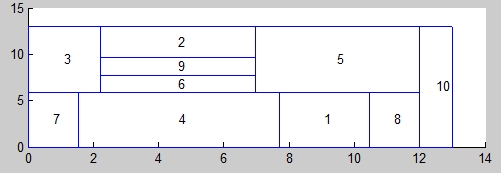}
\end{center}
\caption{Output layout in dynamic condition where U(130,150) is the distribution of costs}
\label{fig:Lay08}
\end{figure}

The five material flows are run and \emph{P-Value} for the experiment is less than 0.0001. In this step, $(\mu+\sigma)$ has been removed and $(\mu+0.25\sigma)$ is added, and in the next iteration, \emph{P-Value} increases to 0.22. Figure \ref{fig:Lay08} shows the best layout according to objective value obtained by simulation which corresponds to material flow $(\mu-\sigma)$.

In this layout rearrange costs have a uniform distribution of U(130 , 150). Because of the fact that rearrange cost is comparatively high, none of departments have changed their location. 
Costs of all layouts are shown in Table \ref{tab:OutSim1}. As shown, the obtained cost of simulation for all material flows is same, which shows the fact that all material flows have converged to the layout without any change.

\begin{table}[H]
\centering
\caption{Simulation results for five layouts with U(130,150) distribution for costs}
\label{tab:OutSim1}
\begin{tabular}{ccc}
\hline
	& $\sum_{ij} f_{ij} d_{ij}$ & 	Simulation \\
\hline
$\mu$ & 351&	{\cellcolor[gray]{0.9}}3518487 \\
$\mu+\sigma$ &300&	{\cellcolor[gray]{0.9}}3518487 \\
$\mu-\sigma$ &{\cellcolor[gray]{0.9}}249&	{\cellcolor[gray]{0.9}}3518487 \\
$\mu+2\sigma$ &274&	{\cellcolor[gray]{0.9}}3518487 \\
$\mu+1.5\sigma$ &287&	{\cellcolor[gray]{0.9}}3518487 \\
\hline
\end{tabular}
\end{table}

On the other hand, rearrangement of departments is related to the comparison of material handling cost in life cycle of products and the rearrangement costs. The algorithm compares the two type of costs and if the rearrangement cost is greater than material handling cost during the product life cycle, then it tries to make the least possible changes.

\begin{figure}
\begin{center}
\includegraphics[width=5in]{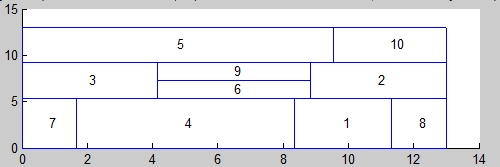}
\end{center}
\caption{The result of change in four restricts}
\label{fig:Lay09}
\end{figure}

It is notable that the maximum ratio of new department is 4 which restricts the changes in the layout in a way that obligates to rearrange some of departments and results in an unfavourable layout in the first run of simulation by best objective function of 4157937 as in Figure \ref{fig:Lay08}. Algorithm could not add new department in the new space of layout without any changes. So, since rearrangement cost and penalty cost do not allow creating an acceptable layout, the maximum ratio for this condition has been increased up to 13, result of which is shown in Figure \ref{fig:Lay09} and Table \ref{tab:OutSim1}.

\subsection{Analysis of results}

The result of our algorithm shows that in a given problem by different amount of rearrangement cost, different material flows result in best layout, i.e. we cannot suggest a given material flow as the best for all problems. This means that in stochastic conditions, considering expected value for material flow distribution as \citet{jithavech2010simulation} suggest or considering a given significance level as \citet{tavakkoli2007design} suggest is not perfect. Although, \citet{tavakkoli2007design} considered normal distribution instead of uniform, but still we cannot suggest a given good material flow that can guarantee good results.%, as the case that we observed in the problem with three different costs. 

Moreover, the sum of the costs of material flow and rearrangement cost of departments severely depends to the amount of material flow. For example, by increasing coefficient $k$ in $\mu +k\sigma$, the material handling cost will increase considerably. If this criterion is considered to decide on the efficiency of layouts, the layout which is the result of the lowest material flow will be selected. We can see, that in result of each scenario in Table \ref{tab:ObjRearrange} the best solution in terms of material handling cost objective is related to the material flow with higher negative coefficient ($k$), but in simulation the efficiency of each layout is calculated considering the random material flows in the product life cycle and it demonstrated that simulation is a justifiable criterion for selecting the best layout.

Finally, we had five initial material flows in the simulation. 
If reasonable computation power is available, selecting more initial material flows could provide quicker high quality results. It is a trade-ff between the computation power and time.

\section{Conclusions}
\label{sec: conclusions}
Unequal Area Facility Layout problem has been one of the most interesting problems in last three decades, since a good layout leads to lower material handling cost in long term. 
In this paper we introduced a DFLP with stochastic material flow, stochastic cost, and stochastic time of change. 

In order to solve this problem, a hybrid simulation with grid search is proposed, in which a discrete event simulation has been developed to realize the market condition. 
In the heart of the algorithm we proposed a genetic algorithm, and to make it quicker we proposed some heuristic solutions to add as initial solutions into the algorithm. As a result, the algorithm automatically starts from the area in which the probability of existence of good solutions is more than an arbitrary areas. 

Our results show that the belief that mean of flows is a good approximate of dynamic material flow in stochastic condition, is not an accurate assumption. In some cases the material flows $\mu –0.625\sigma$ and $\mu-\sigma$ resulted in best solutions, which shows we cannot propose a given material flow to obtain the global optimal layout, and instead our algorithm provides a procedure to achieve good solutions.
Further research can be done by considering the condition of departments with fixed locations, empty spaces in the final layout, departments with constant length and width or corridors in the layout and limiting transportation to them which can make the problem more applicable to the real world conditions.

%\section{Acknowledgments}
%The authors would like to thank ...

\section*{REFERENCES}
\bibliographystyle{ormsv080}
\bibliography{uaflp}

\section{Appendices}

\subsection{Appendix A}

\begin{figure}
\begin{center}
\includegraphics[width=5in]{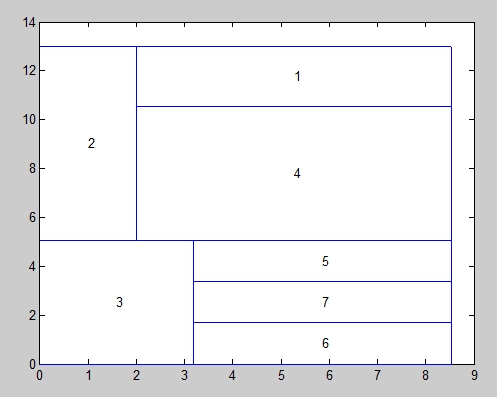}
\end{center}
\caption{Best layout acquired for example 07}
\label{fig:Lay10}
\end{figure}

\begin{figure}
\begin{center}

\end{center}
\centering
\includegraphics[width=5in]{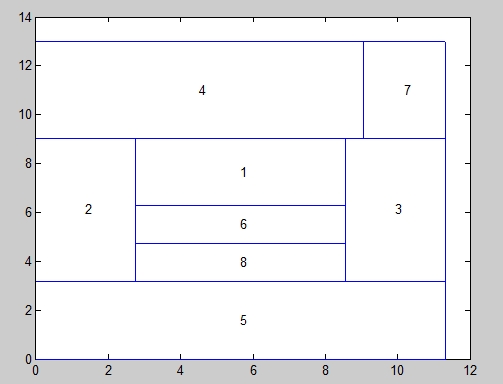}
\caption{Best layout acquired for example 08}
\label{fig:Lay11}
\end{figure}

In order to test the power of our genetic algorithm and specially the result of our heuristic to initiate populations, we tested 5 classic problems from literature. Problem definition is obtained from \citet{komarudin2010applying}. For problems 07 and 08 the best acquired layouts are as Figure \ref{fig:Lay10} and Figure \ref{fig:Lay11} with objective function values of 131.68 and 243.11. This amount for objective function is equal to the best existing amount in literature.

\begin{figure}
\begin{center}
\includegraphics[width=5in]{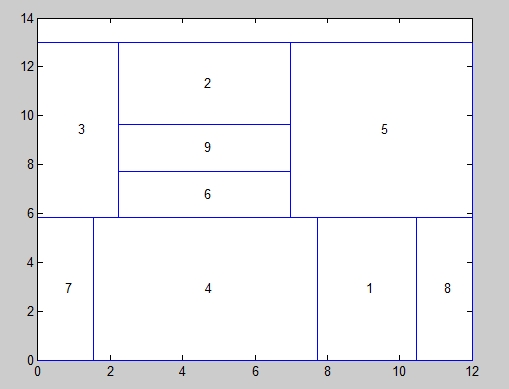}
\end{center}
\caption{Best layout acquired for example 09}
\label{fig:Lay12}
\end{figure}

For problem 09 the best answer acquired is as Figure \ref{fig:Lay12} with objective function of 239.05. This amount for objective function is 1.2\% more than the best existing amount.

For problem AB20 the best answer acquired is as Figure 13 in which the amount of objective function is 5543. This amount for objective function is 9.2 \% more than the best existing amount. 

\begin{figure}
\begin{center}
\includegraphics[width=5in]{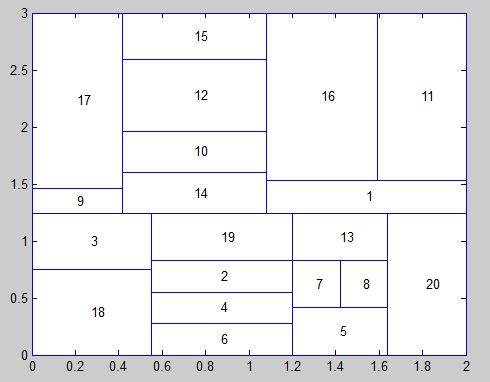}
\end{center}
\caption{Best layout acquired for example AB20}
\label{fig:Lay13}
\end{figure}

\begin{figure}
\begin{center}
\includegraphics[width=5in]{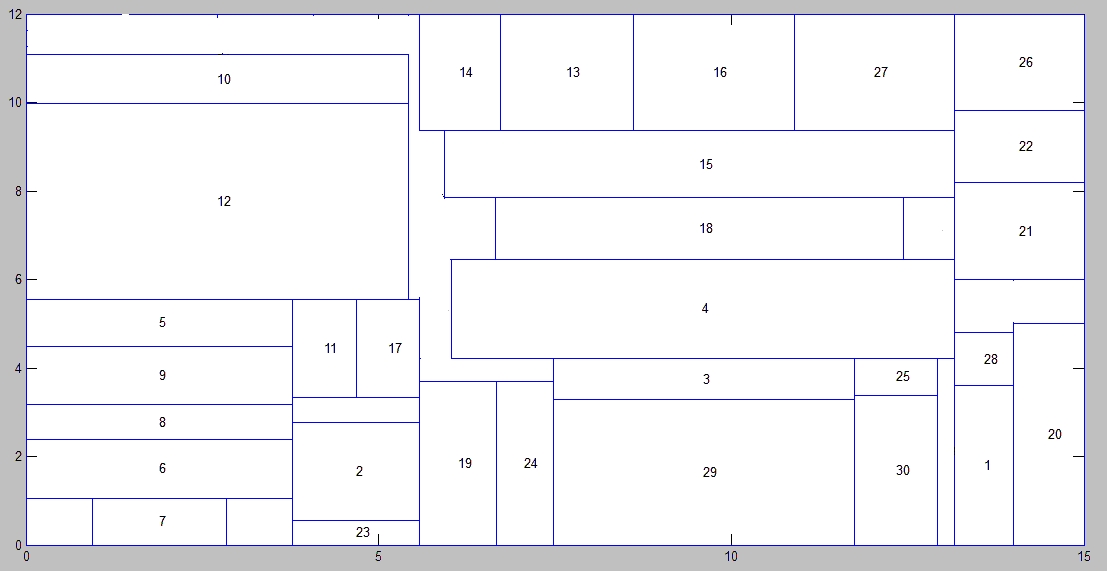}
\end{center}
\caption{Best layout acquired for example SC30}
\label{fig:Lay14}
\end{figure}

For problem SC30 the best answer acquired is as Figure \ref{fig:Lay14} in which the amount of objective function is 4248. This amount for objective function is 14.8 \% more than the best existing amount. 

\end{document}